\documentclass[12pt]{article}
\usepackage{amsmath}
\usepackage{amssymb}
\usepackage{dsfont}
\usepackage{cite}
\newsymbol \blackbox 1004
\newcommand{\eh}{\hfill}\newlength{\sperr}

\newenvironment{proof}{{\settowidth{\sperr}{\bf\rm
Proof}%
\par\addvspace{0.3cm}\noindent\parbox[t]{1.3\sperr}
{\bf\rm P\eh r\eh o\eh o\eh f\eh }%
}}{\nopagebreak\mbox{}
$\blackbox$\par\addvspace{0.3cm}}

\def\nn{\nonumber}
\def\Span{\mathrm{span}\ }

\def\b{\beta}

\def\g{\gamma}
\def\G{\Gamma}

\def\om{\omega}

\def\vp{\varphi}

\def\wh{\widehat}
\def\wt{\widetilde}
\def\ov{\overline}

\def\BC{{\mathbb C}}
\def\BR{{\mathbb R}}

\def\clh{{\mathcal H}}

\def\im{{\rm Im\ }}

\newcommand{\E}{\mathrm{e}}
\newcommand{\I}{\mathrm{i}}

\newtheorem{Pa}{Paper}[section]
\newtheorem{Tm}[Pa]{{\bf Theorem}}
\newtheorem{La}[Pa]{{\bf Lemma}}
\newtheorem{Cy}[Pa]{{\bf Corollary}}
\newtheorem{Rk}[Pa]{{\bf Remark}}

\newtheorem{Dn}[Pa]{{\bf Definition}}

\newtheorem{Pn}[Pa]{{\bf Proposition}}

\newenvironment{dedication}
        {\vspace{1ex}\begin{quotation}\begin{center}\begin{em}}    
        {\par\end{em}\end{center}\end{quotation}}

\title{On accelerants and their analogs, and on the  characterization of the rectangular Weyl functions for Dirac systems
with locally square-integrable potentials on a semi-axis}

\author{Alexander Sakhnovich}

\date{}
\begin{document}
\maketitle

\begin{dedication}
Dedicated  to Heinz Langer      
\end{dedication}

\begin{abstract}    We characterize the set of  rectangular Weyl matrix functions corresponding to 
Dirac systems with locally square-integrable potentials on a semi-axis and
demonstrate  a new way to recover the locally square-integrable potential from the Weyl function.
Important interconnections between our approach and accelerants of convolution operators
are discussed as well.
\end{abstract}



\section{Introduction}\label{Intro}
\setcounter{equation}{0}
\subsection{Some remarks on Dirac system and its accelerant.}
We consider
the self-adjoint Dirac (more precisely, Dirac-type) 
system
\begin{align} &       \label{1.1}
\frac{d}{dx}y(x, z )=\I (z j+jV(x))y(x,
z ) \quad
(x \geq 0),
\end{align} 
where
\begin{align} &   \label{1.2}
j = \left[
\begin{array}{cc}
I_{m_1} & 0 \\ 0 & -I_{m_2}
\end{array}
\right], \hspace{1em} V= \left[\begin{array}{cc}
0&v\\v^{*}&0\end{array}\right],  \quad m_1+m_2=:m,
 \end{align} 
$I_{m_k}$ is the $m_k \times
m_k$ identity
matrix and $v(x)$ is an $m_1 \times m_2$ matrix function.
The matrix functions  $V$ and (sometimes) $v$, are called the { potentials} of Dirac system.
For convenience, we talk about the {\it potential} $v$ in this paper.

Dirac system (also called canonical system, Zakharov-Shabat system, AKNS system and spectral Dirac system)
is a classical object of analysis which has various applications. In particular, it is an auxiliary linear system
for several well-known nonlinear integrable equations (hence, the names Zakharov-Shabat system or AKNS system).
Important interconnections between the Dirac system (spectral Dirac system) and dynamical Dirac system were
studied, for instance, in \cite{ALS-JMP, YaL} (see also the references therein). In addition, Dirac system is a more
general object than the famous Schr\"odinger equation and may be fruitfully used in the study of this equation
(see, e.g., \cite{BolGe, EGNST, GeGo}).

Inverse spectral problem for Dirac system goes back to the seminal paper \cite{Kr55} by M.G. Krein from the year 1955,
where the case of continuous scalar potential $v$ was dealt with.
Numerous works were dedicated to the proofs of the results from \cite{Kr55} and to their generalizations
(see, e.g., \cite{AGKLS, Den, FKRS1,  KrHLa, HLa, LaWo, mykpu, pu, ALS-Dir, RoSa1, RoSa2, SaLoKr} and various references therein).
In particular, the topic is developed further in the interesting recent   Krein-Langer paper \cite{KrHLa}.
An important feature of the Krein's approach is a deep  connection between Dirac systems and convolution operators, which play 
a crucial role in solving inverse problem. The kernel of the corresponding positive convolution (integral) operator  is called {\it accelerant}. 
It proved (see, e.g., Remark 2.3 in \cite{FKS} or Remark 2.46 in \cite{SaSaR})
that accelerant is an analogue of the important $A$-amplitude actively studied and used in the well-known papers
\cite{GeSi, Si}. Moreover, as mentioned in \cite{HLa},
the $A$-amplitude ``is essentially the derivative of Krein's transfer function". Recall that Krein's transfer function appeared in his paper
\cite{Kr53} on Sturm-Liouville equation. Thus,  the accelerant and $A$-amplitude have also Krein's transfer function
as a common forerunner.

Here, we consider the recovery of system from the Weyl function. Such problems form an important area in inverse spectral problems, which
is closely related to the recovery of system from the spectral function.
In the case of Dirac systems \eqref{1.1}  with rectangular matrix functions $v$, one cannot use convolution operators
anymore. However, structured operators  of the form \eqref{r3}  which satisfy operator identities \eqref{r4}, similar to those
for convolution operators, are applied for solving inverse problems \cite{FKRS1, ALS-JST, SaSaR}. (See \cite{SaSaR, SaLopid1, SaL1,
SaL3} and references therein for more details on the method of operator identities.) Somewhat different kinds of structured
operators play the main role in solving inverse problems for skew-self-adjoint Dirac systems, for auxiliary systems in the case of $N$-wave equations
and for certain systems depending rationally on the spectral parameter (see, e.g., \cite{FKRS2, SaA021, ALS111, SaSaR}).
{\it In particular, the matrix function $\Phi_1^{\prime}$ in \eqref{r3} $($or, equivalently, in \eqref{r3!}$)$ may be considered as a direct 
analog of the Krein's accelerant.}

This paper is a continuation of our paper \cite{ALS-JST}, where an inverse problem for Dirac systems
with locally square-integrable potentials on a semi-axis have been solved, and a procedure to recover the $m_1 \times m_2$ potential $v$ from
the Weyl function have been given. The work \cite{ALS-JST} was initiated by a question formulated by F. Gesztesy and generalizes many
results obtained in \cite{FKRS1, ALS-Dir} (for locally bounded potentials) on an essentially wider class of potentials.
In the present paper, we give {\it characterization of Weyl matrix functions} corresponding to such Dirac systems and
demonstrate (see Corollary \ref{CyaltIP}) an alternative way to recover $v$ from the Weyl function.
\subsection{Main notions and results}
We  assume
that  the potential $v$ of Dirac system \eqref{1.1} is  locally square-integrable, that is, square-integrable (or, using an equivalent term, square-summable)
on  the finite intervals $[0, \, l]$. We say that a matrix function is integrable (square-integrable) if its entries are integrable (square-integrable). 
The notation $u(x,z)$ stands for the fundamental solution of \eqref{1.1}
normalized by the condition
\begin{align} &   \label{1.3}
u(0,z)=I_m.
\end{align}
\begin{Dn} \label{defWeyl} A Weyl-Titchmarsh $($or simply Weyl$)$ function of Dirac system  \eqref{1.1} on $[0, \, \infty)$,
where  the potential $v$ is locally integrable, is a holomorphic $m_2\times~m_1$ matrix function $\vp$
which satisfies the inequality
\begin{align} &      \label{2.20}
\int_0^{\infty}
\begin{bmatrix}
I_{m_1} & \vp(z)^*
\end{bmatrix}
u(x,z)^*u(x,z)
\begin{bmatrix}
I_{m_1} \\ \vp(z)
\end{bmatrix}dx< \infty , \quad z\in \BC_+.
\end{align} 
\end{Dn}
Here  $\BC_+$ stands for the open upper half-plane of the complex plane.
\begin{Rk}\label{RkWT} In \cite[Subsection 2.2.1]{SaSaR}, we show that the Weyl function $\vp(z)$ always exists in $\BC_+$ and that it is unique,
holomorphic and contractive $($i.e., $\vp(z)^*\vp(z)\leq I_{m_1})$.
\end{Rk}
Since $\vp$ is contractive, the formula below is well-defined.
\begin{align} \label{5.1}&
\Phi_1\Big(\frac{x}{2}\Big)=\frac{1}{\pi}\E^{x\eta}{\mathrm{l.i.m.}}_{a \to \infty}
\int_{-a}^a\E^{-\I x\zeta}\frac{\vp(\zeta+\I \eta)}{2\I(\zeta +\I \eta)}d\zeta, \quad \eta >0.
\end{align} 
Here l.i.m. stands for the entrywise limit in the norm of  $L^2(0,b)$, 
$\, 0<b \leq \infty$. 
Note that if we put additionally $\Phi_1(x)=0$ for $x<0$, equality  (\ref{5.1})
holds for l.i.m. as the entrywise limit in $L^2(-b,b)$.
The matrix function $\Phi_1$ does not depend on $\eta >0$ (see, e.g., \cite{ALS-JST}).
Moreover, according to the formula (3.12) (which defines $\Phi_1$ in another but equivalent way) and to Proposition 3.1, both from \cite{ALS-JST}, 
the following statement is valid.
\begin{Pn} \label{NesCond} Let $\vp$ be the Weyl function of Dirac system  \eqref{1.1} on $[0, \, \infty)$,
 where the potential $v$ is  locally square-integrable. Then $\Phi_1$ given on $\BR$ by \eqref{5.1} is
 absolutely continuous, $\Phi_1(x)\equiv 0$ for $x\leq 0$, $\Phi_1^{\prime}$ is locally square-integrable on $\BR$, 
 and the operators
\begin{align} &      \label{r3}
S_{\xi}=I-\frac{1}{2}\int_0^{\xi}\int_{|x-t|}^{x+t}\Phi_1^{\prime}\left(\frac{\zeta+x-t}{2}\right)
\Phi_1^{\prime}\left(\frac{\zeta+t-x}{2}\right)^*d\zeta \, \cdot \, dt 
\end{align} 
are positive definite and boundedly invertible in $L^2(0,\xi)$ $(0<\xi <\infty)$.
\end{Pn}
Here (as usual), $\Phi_1^{\prime}:=\frac{d}{d x}\Phi_1$. The results above (on the properties of the Weyl functions) belong  to the direct spectral problem.
The procedure of solving  inverse spectral problem in  \cite{ALS-JST} (see Theorem \ref{TmIP} in the next
section) is based on the construction of  the matrix function $\Phi_1$ and structured operators $S_{\xi}$, which are introduced above.
Below, we consider the functions $\vp(z)$ with the properties as given above and formulate Theorem \ref{MTm} that such functions
are always Weyl functions. In this way, we {\it characterize the set of Weyl functions}.

\begin{Rk}\label{hcontr} It is easy to see that for any $m_2 \times m_1$ matrix function $\vp$,
which is holomorphic and contractive in $\BC_+$, the transformation \eqref{5.1} generates $\Phi_1$
such that $\Phi_1(x) \equiv  0$ for $x\leq 0$. Moreover, this $\Phi_1$ does not depend on $\eta >0$.
\end{Rk}
\begin{Rk}\label{RkOpId} According to \cite[Remark 4.5]{ALS-JST},
the operator $S_{\xi}$ given by \eqref{r3} $($with  an absolutely continuous $m_2 \times m_1$ matrix function  $\Phi_1(x)$ such that $\Phi_1(0)=0$ and
 $\Phi_1^{\prime}(x)$ is square-integrable on $[0,\xi])$ is the unique solution of the operator identity
\begin{align} &      \label{r4}
 A_{\xi}S_{\xi}-S_{\xi}A_{\xi}^*=\I \Pi_{\xi} j \Pi_{\xi}^*,
 \end{align}  
 where $A_{\xi}$ is an integration operator in $L^2(0, \xi)$ multiplied by $-\I$ and  $\Pi_{\xi}$ is a multiplication
 operator $\big(\Pi_{\xi}\in B\big(\BC^m, L^2_{m_2}(0,\xi)\big)\big)$. That is, $A_{\xi}$ and $\Pi_{\xi}$ are given by
 the relations  
 \begin{align} &      \label{r5}
 A_{\xi}=-\I \int_0^x \, \cdot \, dt, \quad \Pi_{\xi}g=\begin{bmatrix} \Phi_1(x) & I_{m_2} \end{bmatrix}g \quad (g\in \BC^m).
 \end{align}  
\end{Rk} 
Changing variables (more precisely, making substitution $\zeta =x+t-2r$) we rewrite $S_{\xi}$ in an equivalent and more
convenient form
\begin{align} &      \label{r3!}
S_{\xi}=I-\int_0^{\xi}s(x,t) \, \cdot \, dt , \quad s(x,t):= \int_{0}^{\min(x,t)}\Phi_1^{\prime}(x-r)
\Phi_1^{\prime}(t-r)^*dr.
\end{align} 

In order to complete a characterization of Weyl functions, we assume that some $m_2 \times m_1$ matrix function
has the properties of the Weyl function described in Remark \ref{RkWT} and Proposition \ref{NesCond}
and prove in this paper the following theorem.
\begin{Tm}\label{MTm}
Let  an $m_2 \times m_1$ matrix function $\vp(z)$ be  holomorphic and contractive  in $\mathbb{C}_+$.
Let $\Phi_1(x)$ given  by \eqref{5.1} be
 absolutely continuous, let $\Phi_1(0)=0$, and let $\Phi_1^{\prime}$ be  square-integrable on 
 all the finite intervals $[0, \xi]$. 
 Assume that the operators $S_{\xi}$, which are expressed via $\Phi_1$ in \eqref{r3!},
are positive definite and boundedly invertible in $L^2(0,\xi)$ $(0<\xi <\infty)$.

Then $\vp$ is the Weyl function
of some Dirac system \eqref{1.1} on $[0, \, \infty)$ such that the potential $v$ of this Dirac system is
locally square-integrable.
\end{Tm}  
In fact, the requirements in Theorem \ref{MTm} slightly differ from the properties in 
Remark \ref{RkWT} and Proposition \ref{NesCond}
but (in view of Remark \ref{hcontr} and equivalence of representations \eqref{r3} and \eqref{r3!}) the
requirements coincide with those properties.

We formulate some results from \cite{ALS-JST} and an auxiliary Proposition \ref{PnV} in the next section
``Preliminaries". The proof of the main Theorem \ref{MTm} and an alternative procedure
of solving inverse problem are contained in Section \ref{Char}.

The notation $\BR$ denotes, in the paper, the real axis, $\BC$ stands for the complex plane and $\BC_+$ stands for the open upper half-plane.
By $B(\clh)$ we denote the set of operators bounded in some Banach space $\clh$, and 
the notation $B(\clh_1,\clh_2)$ stands for the set of bounded operators acting from $\clh_1$ into $\clh_2$.

\section{Preliminaries}\label{Prel}
\setcounter{equation}{0}
The main part of this section, including Theorem \ref{TmIP}, is dedicated to solving inverse problem
and presents some related results from \cite{ALS-JST}.
Consider the fundamental solution $u$ (of Dirac system \eqref{1.1}) at $z=0$ and partition it into block rows:
\begin{align} &      \label{3.1}
\b(x)=\begin{bmatrix}
I_{m_1} & 0
\end{bmatrix}u(x,0), \quad \g(x)=\begin{bmatrix}
0 &I_{m_2}
\end{bmatrix}u(x,0).
\end{align} 
It is shown in \cite{ALS-JST} and easily follows from \eqref{1.1} and \eqref{1.3}  that
\begin{align} &      \label{i5}
\b j \b^*\equiv I_{m_1}, \quad \g j \g^*\equiv -I_{m_2}, \quad \b j \g^*\equiv 0, \quad  \g^{\prime} j \g^*\equiv 0;
\\ &      \label{5.8}
v(x)=\I \b^{\prime}(x)j\g(x)^*.
\end{align} 

Thus, in order to solve inverse problem and recover $v$, it suffices to recover $\b$ and $\g$.
Let $\vp$ be the Weyl function of some Dirac system with a locally square-integrable potential $v$
and let $\Phi_1$ and $S_{\xi}$ be given by \eqref{5.1} and by \eqref{r3}, respectively. Introduce
operators
 \begin{align}
& \label{i10}
\Pi_{\xi}:= \begin{bmatrix}
\Phi_1 & \Phi_2
\end{bmatrix}, \quad 
\Phi_k \in B\big(\mathbb{C}^{m_k}, \, L^2_{m_2}(0, \, \xi)\big);
\\
&      \label{i10'}
\big(\Phi_1 g_1\big)(x)=\Phi_1(x)g_1, \quad
 \Phi_2 g_2=I_{m_2}g_2\equiv g_2,
\end{align}
where $g_k \in \mathbb{C}^{m_k}$ excluding the last $g_2$ in \eqref{i10'}, which stands for 
the natural embedding of $g_2$ into $L^2_{m_2}(0, \, \xi)$.

The Hamiltonian $H=\g\g^*$ may be expressed \cite{ALS-JST} in terms of $\Pi_{\xi}$ and $S_{\xi}$:
 \begin{align} &      \label{i18}
H(\xi)=\g(\xi)\g(\xi)^*=\left(\Pi_{\xi}^*S_{\xi}^{-1}\Pi_{\xi}\right)^{\prime}.
\end{align} 

Now, we may recover $\g$. 
First, for that purpose,  we partition $\g$ into two blocks $\g=\begin{bmatrix}
\g_1 & \g_2
\end{bmatrix}$, where $\g_k$ ($k=1,2$) is an $m_2\times m_k$ matrix function.
Next, we recover  the so called Schur coefficient $\g_2^{-1}\g_1$:                                                                                        
\begin{align} \label{5.2}&
\left(\begin{bmatrix}
0 &I_{m_2}
\end{bmatrix}H\begin{bmatrix}
0 \\ I_{m_2}
\end{bmatrix}\right)^{-1}
\begin{bmatrix}
0 &I_{m_2}
\end{bmatrix}H\begin{bmatrix}
I_{m_1} \\ 0
\end{bmatrix}=(\g_2^*\g_2)^{-1}\g_2^*\g_1=\g_2^{-1}\g_1.
\end{align} 
Here we used the inequality $\det \g_2\not=0$, which follows from the second identity
in \eqref{i5}. The second identity
in \eqref{i5} yields also
$$I_{m_2}-(\g_2^{-1}\g_1)(\g_2^{-1}\g_1)^*=\g_2^{-1}(\g_2^{-1})^*,$$
which implies that the left-hand side of this equality is invertible.
Taking into account  $\det \g_2\not=0$,  we rewrite $\g_1$
 in the form $\, \g_1=\g_2 (\g_2^{-1}\g_1)$,
and the fourth identity in \eqref{i5} we rewrite as 
 $\g_2^{\prime}=\g_1^{\prime}(\g_2^{-1}\g_1)^*$. Therefore, we obtain 
\begin{align} &      \nn
\g_2^{\prime}=(\g_2 (\g_2^{-1}\g_1))^{\prime}(\g_2^{-1}\g_1)^*, \quad {\mathrm{i.e.,}}\\
 &      \label{p18}
\g_2^{\prime}=
\g_2(\g_2^{-1}\g_1)^{\prime}(\g_2^{-1}\g_1)^*
\big(I_{m_2}-(\g_2^{-1}\g_1)(\g_2^{-1}\g_1)^*\big)^{-1},
\end{align} 
and recover $\g_2$ from \eqref{p18} and the initial condition $\g_2(0)=I_{m_2}$.
Finally, we recover $\g_1$ from $\g_2$ and $\g_2^{-1}\g_1$.

In order to recover $\b$ from $\g$, we partition $\b$ into two blocks $\b=\begin{bmatrix}
\b_1 & \b_2
\end{bmatrix}$, where $\b_k$ ($k=1,2$) is an $m_1\times m_k$ matrix function.
We put
\begin{align} \label{5.3}&
\wt \b=\begin{bmatrix}
I_{m_1} & \g_1^*(\g_2^*)^{-1}
\end{bmatrix}.
\end{align} 
Because of  (\ref{i5}) and  (\ref{5.3}), we have $\b j \g^*=\wt \b j \g^*=0$, and so
\begin{align} \label{5.4}&
\b(x)= \b_1(x)\wt \b(x).
\end{align} 
It follows from  (\ref{1.1}) and  (\ref{3.1}) that 
\begin{align} \label{5.4'}&
\b^{\prime}(x)=\I v(x)\g(x),
\end{align} 
which implies
\begin{align} \label{5.5}&
\b^{\prime}j\b^*=0.
\end{align} 
Formula  (\ref{5.4}) and the first relation in  (\ref{i5})  lead us to
\begin{align} \label{5.6}&
 \wt \b j \wt \b^*=\b_1^{-1}(\b_1^*)^{-1}.
\end{align} 
From \eqref{5.4} we also derive that
\begin{align} \nn&
\b^{\prime}j\b^*=\b_1^{\prime}(\wt \b j\wt \b^*)\b_1^*+\b_1(\wt \b^{\prime}j\wt \b^*)\b_1^*.
\end{align}
Taking into account   (\ref{5.5}) and  \eqref{5.6}, we rewrite the
relation above:
\begin{align} \label{5.6'}&
\b_1^{\prime}\b_1^{-1}+\b_1(\wt \b^{\prime}j\wt \b^*)\b_1^*=0.
\end{align}
According to  (\ref{1.3}), \eqref{5.6},  and  (\ref{5.6'}),  $\b_1$ satisfies the first order differential equation
(and initial condition):
\begin{align} \label{5.7}&
\b_1^{\prime}=-\b_1(\wt \b^{\prime}j\wt \b^*)( \wt \b j \wt \b^*)^{-1},  \quad   \b_1(0)=I_{m_1}.
\end{align} 
Thus, $\b_1$ and $\b$ are successively recovered from $\g$.
The potential $v$ is recovered from $\b$ and $\g$ via  \eqref{5.8}.
We obtain the following theorem.
 \begin{Tm} \label{TmIP} Let $\vp$ be the Weyl function of Dirac system  \eqref{1.1} on $[0, \, \infty)$,
 where the potential $v$ is  locally square-integrable.
 Then $v$ can be uniquely recovered from $\vp$ via the formula
\eqref{5.8}.
Here, $\b$ is recovered from $\g$ using  \eqref{5.3},  \eqref{5.4} and  \eqref{5.7}$;$ $\g$ is recovered
from the Hamiltonian $H$ using  \eqref{5.2} and   \eqref{p18}$;$ the Hamiltonian is given by
 \eqref{i18}, $\Pi_{\xi}$  from  \eqref{i18} is expressed via $\Phi_1(x)$ in formulas \eqref{i10} and
 \eqref{i10'}, and  $S_{\xi}$ is  expressed via $\Phi_1(x)$ in \eqref{r3}. Finally, $\Phi_1(x)$ is recovered from $\vp$ using  \eqref{5.1}.
 \end{Tm}
 In a way, which is similar to the recovery of $\b$ from $\g$ (in Theorem \ref{TmIP}), we recover $\g$ from $\b$
 in the next proposition.  This proposition (see below) is a simple modification of  \cite[Proposition 2.53]{SaSaR} for the
 case of an absolutely continuous matrix function $\b$ with the proof of \cite[Proposition 2.53]{SaSaR}  remaining valid
 for our case. (We note that we use the notation $\wh \g$ instead of $\g$ in our proposition.)
 \begin{Pn} \label{PnV} Let a given $m_1 \times m$ matrix function $\b(x)$ $($$0 \leq x \leq \xi$$)$
be absolutely continuous and satisfy relations
\begin{align} &      \label{b1}
\b(0)=\begin{bmatrix}
I_{m_1} & 0
\end{bmatrix}, \qquad \b^{\prime}j\b^* \equiv 0.
\end{align}  
Then there is a unique  $m_2 \times m$ matrix function $\wh \g$,
which is absolutely continuous  and satisfies relations
\begin{align} &      \label{b2}
\wh \g(0)=\begin{bmatrix}
0 & I_{m_2}
\end{bmatrix}, \quad \wh \g^{\prime}j \wh \g^* \equiv 0, \quad  \wh \g j\b^*\equiv 0.
\end{align}  
This $\wh \g$ is given by the formula
\begin{align} &      \label{b3}
\wh \g=\wh \g_2\wt \g, \quad \wt \g:=\begin{bmatrix}
\wt \g_1 & I_{m_2}
\end{bmatrix}, \quad \wt \g_1:=\b_2^*(\b_1^*)^{-1},
\end{align}  
where  $\wh \g_2$ is recovered via  the differential system and initial condition below$:$
\begin{align} &      \label{b4}
\wh \g_2^{\prime}=\wh \g_2\wt \g_1^{\prime} \wt \g_1^*\big(I_{m_2}- \wt \g_1  \wt \g_1^*\big)^{-1}, \quad \wh \g_2(0)=I_{m_2}.
\end{align} 
Moreover, the procedure above is well-defined since 
\begin{align} &      \label{b5}
 \det \b_1(x)\not=0, \quad \det\big(I_{m_2}- \wt \g_1(x)  \wt \g_1(x)^*\big)\not=0.
\end{align} 
 \end{Pn}
\section{Characterization of Weyl functions}\label{Char}
\setcounter{equation}{0}
In order to prove Theorem \ref{MTm} we need the following simple lemma.
\begin{La} \label{La3.1} Let $\Phi_1^{\prime}$ be square-integrable on $[0, \xi]$ and let the operator $S_{\xi}$ of the
form \eqref{r3!} be positive definite in $L^2_{m_2}(0,\xi)$ and have a bounded inverse.
Then $S_{\xi}^{-1}$ admits a unique factorization
\begin{align} &      \label{c1}
S_{\xi}^{-1}=E_{\Phi,\xi}^*E_{\Phi,\xi}, \quad E_{\Phi,\xi}=I+\int_0^x E_{\Phi}(x,t)\, \cdot \, dt \in
B\big(L^2_{m_2}(0,\xi)\big),
\end{align}  
where $E_{\Phi}(x,t)$ is continuous $($with respect to $x$ and $t)$.
\end{La}
\begin{proof}. First, let us show that the kernel $s(x,t)$ given by \eqref{r3!} is continuous.
Since $\Phi_1^{\prime}$ is square-integrable, its entries may be approximated
in $L^2(0,\xi)$ by some continuous matrix functions $\Psi_{\Delta}(x)$ such that the norms
of all the differences between the corresponding entries of 
$\Phi_1^{\prime}$ and $\Psi_{\Delta}$ are less than $\Delta$.
Then the differences between  $s(x,t)$ and
continuous matrix functions $$s_{\Delta}(x,t):=\int_0^{\min(x,t)}\Psi_{\Delta}(x-r)\Psi_{\Delta}(t-r)^*dr$$
have the form
 \begin{align}     \label{c2}
 s(x,t)-s_{\Delta}(x,t)=\int_0^{\min(x,t)}\big(&(\Phi_1^{\prime}(x-r)-\Psi_{\Delta}(x-r))\Phi_1^{\prime}(t-r)^*
\\ &  \nn
+\Psi_{\Delta}(x-r)
(\Phi_1^{\prime}(t-r)-\Phi_{\Delta}(t-r))^*\big)dr,
\end{align} 
and these differences may be made uniformly (with respect to $x$ and $t$) sufficiently small. It follows that
$s(x,t)$ is, indeed, continuous.

Now, the statement of the lemma follows from the results in \cite[pp. 184-186]{GoKr}
(see also a shorter factorization Corollary 1.39 in \cite{SaSaR}).
\end{proof}
\begin{Rk} \label{Rk3.2} Similar to the considerations from \cite[p. 34]{SaSaR}, we show that the kernel of the integral operator $E_{\Phi,\xi}$ does not depend on $\xi$.  
Indeed, assume that
$\Phi_1$ is given on $[0,\, \ell]$ $(\ell >\xi)$ and satisfies the conditions of Lemma~\ref{La3.1} on $[0,\ell]$.
Then, $S_{\xi}$  satisfies the conditions of Lemma~\ref{La3.1} as well, and $S_{\xi}^{-1}$ and $S_{\ell}^{-1}$ admit unique
factorizations \eqref{c1} and
\begin{align} &      \label{f0}
S_{\ell}^{-1}=E^*E, \quad E=I+\int_0^x E(x,t)\, \cdot \, dt \in
B\big(L^2_{m_2}(0,\ell)\big),
\end{align}  
respectively. In view of \eqref{f0},
for the projector $P_{\xi}\in B\big(L^2_{m_2}(0,\ell), \, L^2_{m_2}(0,\xi)\big)$, such that
\begin{align} &      \label{f1}
 \big(P_{\xi}f\big)(x)=f(x) \quad (0<x<\xi, \quad f\in L^2_{m_2}(0,\ell)),
\end{align}  
we have
\begin{align} &      \label{f2}
S_{\xi}=P_{\xi}S_{\ell}P_{\xi}^* =P_{\xi}E^{-1}(E^*)^{-1} P_{\xi}^*.
\end{align}  
Since $E^{\pm 1}$ are lower triangular operators and $P_{\xi}$ has the form \eqref{f1}, the 
following equalities are valid: 
\begin{align} &      \label{f3}
(E^*)^{-1} P_{\xi}^*=P_{\xi}^*P_{\xi}(E^*)^{-1} P_{\xi}^*, \quad P_{\xi}(E^*)^{-1} P_{\xi}^*=(P_{\xi}E^* P_{\xi}^*)^{-1}.
\end{align}  
Using \eqref{f2} and \eqref{f3}, we obtain
\begin{align} &      \label{f4}
S_{\xi}^{-1}=(P_{\xi}E^* P_{\xi}^*) (P_{\xi}E P_{\xi}^*).
\end{align}  
Now, it is immediate from the uniqueness of the factorization \eqref{c1} that
$E_{\Phi, \xi}=P_{\xi}E P_{\xi}^*$ or, equivalently, $E_{\Phi}(x,t)=E(x,t)$ for
$x \leq \xi$. Thus, we see that $E_{\Phi}(x,t)$ in \eqref{c1} does not depend on $\xi$.
\end{Rk}
Further, we modify the proof of \cite[Theorem 5.2]{FKRS1}  (see also \cite[Theorem 2.54]{SaSaR}) and
 introduce the matrix functions $\b_{\Phi}$ and $\g_{\Phi}$ by  $\b_{\Phi}(0):=\b_{\Phi}(+0)$,
\begin{align} &      \label{c3}
\b_{\Phi}(x):=\begin{bmatrix}I_{m_1} &0 \end{bmatrix}
+\int_0^x\Big(S_x^{-1}\Phi_1^{\prime}\Big)(t)^*\begin{bmatrix}\Phi_1(t) & I_{m_2} \end{bmatrix}dt \quad (x>0);
\\ &      \label{c4} \g_{\Phi}(x):=\begin{bmatrix}\Phi_1(x) & I_{m_2} \end{bmatrix}
+\int_0^xE_{\Phi}(x,t)\begin{bmatrix}\Phi_1(t) & I_{m_2} \end{bmatrix}dt \quad (x \geq 0),
\end{align}  
where $S_x^{-1}$ is applied to $\Phi_1^{\prime}$ in \eqref{c3}  columnwise.
To proceed with our proof of Theorem \ref{MTm}   we prove first the lemma below. 
\begin{La} \label{LaBeGa} 
Let an $m_2 \times m_1$ matrix function $\Phi_1(x)$  be
 absolutely continuous on $[0, \infty)$, let $\Phi_1(0)=0$, and let $\Phi_1^{\prime}$ be  locally square-integrable. 
Assume that the
operators $S_{\xi}$, which are expressed via $\Phi_1$
in  \eqref{r3!}, 
are boundedly invertible for all $\, 0<\xi < \infty$.
Then,  $\b_{\Phi}$ is  
absolutely continuous on $[0, \infty)$, $\g_{\Phi}$ is continuous, and $\b_{\Phi}$ and $\g_{\Phi}$ satisfy the conditions
\begin{align} &      \label{c5}
\b_{\Phi}(0)=\begin{bmatrix}I_{m_1} &0 \end{bmatrix}, \quad
\b_{\Phi}^{\prime}j\b_{\Phi}^*\equiv 0; 
\\ &      \label{c6} 
\g_{\Phi}(0)=\begin{bmatrix}0& I_{m_2} \end{bmatrix}, \quad
 \g_{\Phi}j\b_{\Phi}^*\equiv 0.
\end{align}  
\end{La}
\begin{proof}.
Factorizing $S_x^{-1}$,  we rewrite
\eqref{c3} in the form
  \begin{align}       \nn
 \b_{\Phi}(x) &
 =\begin{bmatrix}I_{m_1} &0 \end{bmatrix}+
 \int_0^x \big(E_{\Phi,x}\Phi_1^{\prime}\big)(t)^*\Big(E_{\Phi,x}\begin{bmatrix}\Phi_1 & I_{m_2} \end{bmatrix}\Big)(t)dt
\\  \label{c3'} &
=\begin{bmatrix}I_{m_1} &0 \end{bmatrix}+
 \int_0^x \big(E_{\Phi, \xi}\Phi_1^{\prime}\big)(t)^*\Big(E_{\Phi,\xi}\begin{bmatrix}\Phi_1 & I_{m_2} \end{bmatrix}\Big)(t)dt
\end{align}
for $\xi  \geq x$. (Here, we took into account Remark \ref{Rk3.2} in order to substitute $E_{\Phi, \xi}$ instead of $E_{\Phi, x}$.)
Clearly,  (\ref{c4}) is equivalent to the equalities
\begin{align} &      \label{c7}
 \g_{\Phi}(x)=\Big(E_{\Phi,\xi}\begin{bmatrix}\Phi_1 & I_{m_2} \end{bmatrix}\Big)(x),
 \quad 0 \leq x \leq \xi \quad ({\mathrm{for}} \,\, {\mathrm{all}} \,\, \xi <\infty),
\end{align}  
where $E_{\Phi,\xi}$ is given in \eqref{c1}. 
We note that, according to Remark \ref{Rk3.2} and to the conditions of lemma, 
the matrix functions $\b_{\Phi}$ and $\g_{\Phi}$ given by \eqref{c3} and \eqref{c4}
are well-defined for all $x \geq 0$. 
On the other hand, in order to prove \eqref{c5} and \eqref{c6} it suffices
to  prove that \eqref{c5} and \eqref{c6} hold  for   all $x<\xi$ with any arbitrary fixed $\xi >0$, and so we fix some arbitrary $\xi$.

The first equalities in  (\ref{c5}) and   (\ref{c6})
are immediate from   (\ref{c3}) and   (\ref{c4}), respectively.
Next, we multiply both sides of \eqref{r4} by $E_{\Phi,\xi}$ from the left and by $E_{\Phi,\xi}^*$ from the right.
Taking into account  \eqref{c1}, we obtain 
$$E_{\Phi,\xi}A_{\xi}E_{\Phi,\xi}^{-1}-\big(E_{\Phi,\xi}^{-1}\big)^*A_{\xi}^*E_{\Phi,\xi}^*=\I E_{\Phi,\xi}\Pi_{\xi}j\Pi_{\xi}^*E_{\Phi,\xi}^*.$$
Using the expression for $\Pi_{\xi}$ in \eqref{r5} and formula \eqref{c7}, we rewrite the equality above in the form
\begin{align} &\nn
E_{\Phi,\xi}A_{\xi}E_{\Phi,\xi}^{-1}-\big(E_{\Phi,\xi}^{-1}\big)^*A_{\xi}^*E_{\Phi,\xi}^*=\I\g_{\Phi}(x)j\int_0^{\xi}\g_{\Phi}(t)^*\,\cdot \,dt,
\end{align}
which (taking into account that $E_{\Phi,\xi}A_{\xi}E_{\Phi,\xi}^{-1}$ is a triangular operator) yields
 \begin{align} &      \label{c8}
E_{\Phi,\xi}A_{\xi}E_{\Phi,\xi}^{-1} =\I\g_{\Phi}(x)j\int_0^x\g_{\Phi}(t)^*\,\cdot \,dt.
 \end{align}  
Let us partition $\g_{\Phi}$ into the blocks: 
$$\g_{\Phi}=\begin{bmatrix} (\g_{\Phi})_1 & (\g_{\Phi})_2 \end{bmatrix}, \quad {\mathrm{where}} \quad
(\g_{\Phi})_1(x)=\big(E_{\Phi,\xi}\Phi_1\big)(x) \quad {\mathrm{for}} \quad x \leq \xi.$$ 
Using expression  \eqref{r5} for $A_{\xi}$ and the equality $\Phi_1(0)=0$, we see that
 \begin{align} &      \label{c9}
\big(E_{\Phi,\xi}A_{\xi}E_{\Phi,\xi}^{-1}\big) E_{\Phi,\xi}\Phi_1^{\prime}=-\I E_{\Phi,\xi}\Phi_1=-\I (\g_{\Phi})_1.
 \end{align}  
Relations \eqref{c7}, \eqref{c8} and \eqref{c9} imply the equality
 \begin{align}       \nn
 (\g_{\Phi})_1(x)& =-\g_{\Phi}(x)j\int_0^x\g_{\Phi}(t)^*\big(E_{\Phi,\xi}\Phi_1^{\prime}\big)(t)dt
 \\  \label{c10} &
 =-\g_{\Phi}(x)j
 \int_0^x\Big(E_{\Phi,\xi}\begin{bmatrix}\Phi_1 & I_{m_2} \end{bmatrix}\Big)(t)^*\big(E_{\Phi,\xi}\Phi_1^{\prime}\big)(t)dt.
 \end{align}  
From \eqref{c3'} and \eqref{c10} we derive
 \begin{align}    \label{c11} &
 (\g_{\Phi})_1(x)
 =-\g_{\Phi}(x)j
 \big(\b_{\Phi}(x)-\begin{bmatrix}I_{m_1} &0 \end{bmatrix}\big)^*.
 \end{align}  
It is immediate that \eqref{c11} is equivalent to
  \begin{align}      &  \label{c12}
 \g_{\Phi}(x)j\b_{\Phi}(x)^* \equiv 0.
 \end{align} 
According to relations \eqref{c3'} and \eqref{c7}, $\b_{\Phi}$ is absolutely continuous and 
almost everywhere on $[0, \xi]$ we have
\begin{align} &      \label{c13}
\b_{\Phi}^{\prime}(x)=
\Big(E_{\Phi,\xi}\Phi_1^{\prime}\Big)(x)^*\g_{\Phi}(x).
 \end{align} 
The second equality in \eqref{c5} easily follows from \eqref{c12} and \eqref{c13}.
Thus,  \eqref{c5} is proved. We derived the first equality in \eqref{c6} at the beginning of lemma's proof whereas  the second equality 
in \eqref{c6} coincides with \eqref{c12}, that is, \eqref{c6} holds as well.
\end{proof}
Now, we can prove our main theorem.

{\it Proof of Theorem \ref{MTm}.} 
Step 1. Let us consider $\b_{\Phi}$ and $\g_{\Phi}$ (constructed via $\Phi_1$) in greater detail
than in Lemma \ref{LaBeGa}.
We note that, according to \eqref{c13}, $\b_{\Phi}^{\prime}$ is locally square-integrable
and so $\wh \g$ corresponding to $\b_{\Phi}$ and given by the formulas \eqref{b3} and \eqref{b4}
in Proposition \ref{PnV} has a locally square-integrable derivative $\wh \g^{\prime}$.
From Proposition \ref{PnV} we see that $\wh \g$ satisfies relations
\begin{align} &      \label{c14}
\wh \g(0)=\begin{bmatrix}
0 & I_{m_2}
\end{bmatrix}, \quad \wh \g^{\prime}j \wh \g^* \equiv 0, \quad  \wh \g j\b_{\Phi}^*\equiv 0.
\end{align}  
First, we  show that 
$\g_{\Phi}$ coincides with $\wh \g$ and so $\g_{\Phi}$ is absolutely continuous,
$\g_{\Phi}^{\prime}$  is  locally square-integrable and   the  equality
\begin{align} &      \label{c15}
\g_{\Phi}^{\prime}j\g_{\Phi}^* = 0
 \end{align} 
is valid (in addition to equalities \eqref{c6} which are already proved).

Since $E_{\Phi}(x,t)$ is continuous and does not depend on $\xi$, the resolvent kernel $\G_{\Phi}$ of 
$$E_{\Phi,\xi}^{-1}=I+\int_0^x\G_{\Phi}(x,t)\cdot dt$$ 
is continuous and does not depend on $\xi$ as well.  
We rewrite  \eqref{c8} in the form of an equality for kernels:
\begin{align}      \nn
I_{m_2}+\int_t^x\big(E_{\Phi}(x,r)+\G_{\Phi}(r,t)\big)dr+
\int_t^x\int_{\xi}^x & E_{\Phi}(x,r)dr\G_{\Phi}(\xi,t)d\xi
\\ &  \label{c}
=-\g_{\Phi}(x)j\g_{\Phi}(t)^*.
 \end{align}  
 In particular, formula  (\ref{c}) for the case $x=t$ implies that
 \begin{align}      &  \label{c0}
 \g_{\Phi}(x)j\g_{\Phi}(x)^* \equiv - I_{m_2}.
 \end{align}  
Formulas  \eqref{c5} and  (\ref{c14})  yield useful relations
  \begin{align} &      \label{c16}
 \b_{\Phi}j\b_{\Phi}^*\equiv I_{m_1}, \quad 
\wh \g j \wh \g^*\equiv - I_{m_2}.
 \end{align}  
 From \eqref{c6}, \eqref{c0}, \eqref{c16} and the last equality in \eqref{c14} we derive that
 \begin{align} &      \label{c22}
 \g_{\Phi}(x)=\om(x) \wh \g(x), \quad \om(x)^*= \om(x)^{-1}
 \end{align}  
for some $m_2\times m_2$ matrix function $\om$.  Similar to the proof of Lemma \ref{LaBeGa} we fix some arbitrary $\xi >0$. We
 show that  $\om(x) \equiv I_{m_2}$ on $[0,\xi]$.

Indeed, in view of the first two equalities in \eqref{c14} and the second equality in \eqref{c16},
 using \cite[Proposition 2.1]{ALS-JST} and the proof of  \cite[Proposition 3.1]{ALS-JST},
 we see that there is an operator 
 \begin{align} &      \label{c20}
\wh E=I+\int_0^x N(x,t)\, \cdot \, dt \in
B\big(L^2_{m_2}(0,\xi)\big),
\end{align}  
 such that 
 \begin{align} &      \label{c21}
\wh E A_{\xi}=\I \wh \g(x)j\int_0^x\wh \g(t)^* \, \cdot \, dt\, \wh E, \quad \wh \g_2=\wh E I_{m_2}.
 \end{align}  
 Moreover, $N(x,t)$ is a Hilbert-Schmidt kernel and the operators $\wh E^{\pm 1}$ map differentiable
 functions with a square-integrable derivative into differentiable
 functions with a square-integrable derivative. (Above, we  repeat, after some renaming, the statement of 
  \cite[Proposition 3.1]{ALS-JST}, since its proof does not depend really on the existence of the Dirac system
  but follows from  \eqref{c14} and  \eqref{c16}.)
 According to  (\ref{c22}) and  (\ref{c21}), for
  \begin{align} &      \label{c23}
\wh E_{\Phi, \xi}:=\om(x) \wh E
 \end{align}  
  we have
 \begin{align} &      \label{c24}
\wh E_{\Phi, \xi}A_{\xi}=\I \g_{\Phi}(x) j\int_0^x \g_{\Phi}(t)^* \, \cdot \, dt\, \wh E_{\Phi, \xi}, \quad  
(\g_{\Phi})_2=\wh E_{\Phi,\xi}I_{m_2}.
 \end{align}   
On the other hand, formulas  (\ref{c7}) and  (\ref{c8}) lead us to
 \begin{align} &      \label{c25}
 (\g_{\Phi})_2=E_{\Phi,\xi}I_{m_2}, \quad
 E_{\Phi,\xi}A_{\xi}=\I \g_{\Phi}(x)j\int_0^x \g_{\Phi}(t)^* \, \cdot \, dt \, E_{\Phi,\xi}.
\end{align}  
It is easy to see that 
$$\ov{\Span}\Big(\bigcup_{i=0}^{\infty}\im \big(A_{\xi}^i I_{m_2}\big)\Big)=
L^2_{m_2}(0,\xi),$$
where $\im$ stands for image and $\ov{\Span}$ stands for the closed linear span.
Hence, equalities   (\ref{c24}) and  (\ref{c25})
imply that $\wh E_{\Phi,\xi}=E_{\Phi,\xi}$.  Therefore, comparing  the representation of $E_{\Phi,\xi}$ in (\ref{c1})
and formulas (\ref{c20})
and  (\ref{c23})  for $\wh E_{\Phi,\xi}$ we obtain $\om(x) \equiv I_{m_2}$. In other words we have
$\g_{\Phi}\equiv \wh \g$. Thus, $\g_{\Phi}$ is absolutely continuous,
$\g_{\Phi}^{\prime}$  is  locally square-integrable and   \eqref{c15} holds.

Step 2.  From the second equalities in \eqref{c5} and \eqref{c6} and from \eqref{c15}, we derive
\begin{align}   &      \label{c17}
 u_{\Phi}^{\prime}j   u_{\Phi}^* j=\I j
\begin{bmatrix}0&  v \\  v^* &0\end{bmatrix} \quad
{\mathrm{for}} \quad  u_{\Phi}(x):= \begin{bmatrix}\b_{\Phi}(x) \\  \g_{\Phi}(x) \end{bmatrix},
 \quad
  v:=\I \b_{\Phi}^{\prime}j \g_{\Phi}^*.
  \end{align}  
Moreover, relations \eqref{c6}, \eqref{c0} and the first equalities in \eqref{c5},  \eqref{c16}  imply
that 
\begin{align}   &      \label{c18}
 u_{\Phi}(0)=I_m, \qquad u_{\Phi} j   u_{\Phi}^* \equiv  j.
  \end{align}  
 In particular, since   $u_{\Phi}$ is $j$-unitary, we rewrite \eqref{c17} in the form
  \begin{align}   &      \label{c19}
u_{\Phi}^{\prime}=\I j
\begin{bmatrix}0&  v \\  v^* &0\end{bmatrix}  u_{\Phi}.
 \end{align}  
 Hence, $u_{\Phi}(x)$ is the normalized fundamental solution (at $z=0$) of Dirac system with the potential
 $v$ of the form \eqref{c17}. We note that $\wh E=E_{\Phi,\xi}$ is constructed for this system in precisely the same way
 as the operator $E$ is constructed in \cite{ALS-JST}. Recall that $\Phi_1$ is introduced as $E^{-1}\g_1$ in \cite[f-la (3.12)]{ALS-JST}, and note that, in the present paper, 
 the equality $\Phi_1=E_{\Phi,\xi}^{-1}(\g_{\Phi})_1$
 follows from \eqref{c7}. Thus, we may apply the results from \cite{ALS-JST} to our system (with $v=\I \b_{\Phi}^{\prime}j \g_{\Phi}^*$) and to our $\Phi_1$.
 Denoting the Weyl function of this system by $\vp_{\Phi}$, we obtain from \cite[f-la (4.18)]{ALS-JST} the representation
 \begin{align} \label{c26}&
\Phi_1\Big(\frac{x}{2}\Big)=\frac{1}{\pi}\E^{x\eta}{\mathrm{l.i.m.}}_{a \to \infty}
\int_{-a}^a\E^{-\I x\zeta}\frac{\vp_{\Phi}(\zeta+\I \eta)}{2\I(\zeta +\I \eta)}d\zeta, \quad \eta >0.
\end{align} 
 (The existence of the Weyl function is stated in  \cite[Proposition 1.3]{ALS-JST} and is discussed in detail in \cite[Subsection 2.2.1]{SaSaR}.)
 
 The initial function $\vp$ from the main Theorem \ref{MTm}, which we are proving here, generates the same $\Phi_1$ via formula \eqref{5.1}. Subtracting both sides
 of \eqref{c26} from the corresponding sides of \eqref{5.1}, we derive:
  \begin{align} \label{c27}&
{\mathrm{l.i.m.}}_{a \to \infty}
\int_{-a}^a\E^{-\I x\zeta}\frac{\vp(\zeta+\I \eta)-\vp_{\Phi}(\zeta+\I \eta)}{2\I(\zeta +\I \eta)}d\zeta, \quad \eta >0,
\end{align} 
where (taking into account Remark \ref{hcontr}) l.i.m. stands for the entrywise limit in the norms of  $L^2(-b,b)$
($\, 0< b \leq \infty$). Therefore, we see that $\vp_W=\vp$, that is, 
$\vp$ is the Weyl function of the Dirac system  with the potential $v=\I \b_{\Phi}^{\prime}j \g_{\Phi}^*$. $\Box$

The proof of Theorem \ref{MTm} yields an alternative way of recovering $\b$ and $\g$ while solving the
inverse problem, which is different from the one presented in \cite{ALS-JST} (and in Theorem \ref{TmIP} in ``Preliminaries"). 
\begin{Cy}\label{CyaltIP}
Let $\vp$ be the Weyl function of Dirac system  \eqref{1.1} on $[0, \, \infty)$,
 where the potential $v$ is  locally square-integrable.
 
 Then $v$ can be uniquely recovered from $\vp$  using the following procedure.
 First, $\Phi_1(x)$ is recovered from $\vp$ using  \eqref{5.1}. Next, the operators $S_{\xi}$
 are expressed via $\Phi_1(x)$ in formula \eqref{r3} or, equivalently, in \eqref{r3!}. Finally, we set
\begin{align} &      \label{c28}
\b(x)=\begin{bmatrix}I_{m_1} &0 \end{bmatrix}
+\int_0^x\Big(S_x^{-1}\Phi_1^{\prime}\Big)(t)^*\begin{bmatrix}\Phi_1(t) & I_{m_2} \end{bmatrix}dt,
\end{align}
and recover $\g$ from $\b$ using the procedure  from Proposition \ref{PnV} for the recovery of $\wh \g$ $($and putting $\g=\wh \g)$.
The potential $v$ is expressed $($in \eqref{5.8}$)$ via $\b$ and $\g$, namely, $v=\I \b^{\prime}j\g^*$.

\end{Cy}

\bigskip

\noindent{\bf Acknowledgments.}
 {This research   was supported by the
Austrian Science Fund (FWF) under Grant  No. P29177.}

\begin{flushright}

A.L. Sakhnovich,\\
Fakult\"at f\"ur Mathematik, Universit\"at Wien, \\
Oskar-Morgenstern-Platz 1, A-1090 Vienna, Austria\\
e-mail: oleksandr.sakhnovych@univie.ac.at
\end{flushright}


\end{document}